\documentclass{amsart}

\usepackage{amssymb}
\usepackage{hyperref}

\newtheorem{thm}{Theorem}[section]
\newtheorem{cor}[thm]{Corollary}

\newtheorem{lem}[thm]{Lemma}
\theoremstyle{definition}

\newtheorem{ex}[thm]{Example}
\numberwithin{equation}{section}

\begin{document}  

\title{Desingularizations of Some Weighted Projective Planes}
\author{Jeremiah M. Kermes}
\address{5022 Holly Ridge Dr.\\ Raleigh, NC 27612}
\email{cermz27@bellsouth.net}

\begin{abstract}
In this paper we discuss the desingularization algorithm for a toric surface.  In particular, we construct an iterable method of determining the Hirzebruch-Jung continued fraction decomposition.  These results are then applied to weighted projective planes with at least one tivial weight, ${\mathbb P}(1,m,n)$.  The paper concludes with the development of a computer program that computes this continued fraction decomposition.  
\end{abstract}

\maketitle

\section{The Desingularization Algorithm} 
\label{sec:algorithm} 

In this section we review an algorithm in \cite{Oda} for desingularizing a toric surface.  To begin let $\sigma=\left<v_1,v_2\right>$ be a cone in ${\mathbb Z}^2$.  It is known from \cite{Fult} that $U_\sigma$ is non-singular if and only if $\det|v_1\: v_2|=\pm1$.  

In case $\sigma$ yields a singular surface one can desingularize it by subdividing $\sigma$ into smaller cones.  These subdivisions occur by taking rays through a set of lattice points $\left\{l_0,\ldots,l_{s+1}\right\}$ with $l_0=v_1$ and $l_{s+1}=v_2$.  

To begin the algorithm one must obtain 3 pieces of data from $v_1$ and $v_2$.  The first is a lattice point $n_1\in{\mathbb Z}^2$, while the second two are relatively prime integers, $0\leq p_0<q_0$.  The lattice point $n_1$ and $l_0=v_1$ must form a basis for ${\mathbb Z}^2$, and their cone must contain $v_2$.  The integers $p_0$ and $q_0$ must also satisfy \begin{equation}v_2=p_0l_0+q_0n_1.\end{equation}  

With this data in hand, one can construct integers $b_1,\ldots,b_s$ with each $b_j\geq2$.  Begin by forming the fraction $\beta_0=\frac{q_0}{q_0-p_0}>1$.  From here one needs to obtain the Hirzebruch-Jung continued fraction expansion of $\beta_0$.  This means expressing the fraction as \begin{equation} \beta_0=b_1-\frac{1}{\displaystyle b_2-\frac{1}{\displaystyle \ddots -\frac{1}{\displaystyle b_s}}}\end{equation} or $\beta_0=[[b_1,\ldots,b_s]]$ for short.  

Once these integers have been obtained, one constructs both the $l_j$'s for $0\leq j\leq s+1$ and a set of lattice points $n_1,\ldots,n_{s+1}$ such that $l_j$ and $n_{j+1}$ form a ${\mathbb Z}$-basis for ${\mathbb N}$.  The process is inductive, and begins with $l_0=v_1$ and $n_1$ as above.  Then one lets \begin{equation}\begin{gathered} l_{j+1}=l_j+n_{j+1} \\ n_{j+1}=(b_j-2)l_{j-1}+(b_j-1)n_j.\end{gathered}\end{equation}  The rays through $l_j$ for $0\leq j\leq s+1$ form $\Delta(1)$ for the minimal desingularization of $U_\sigma$.  The maximal cones are formed by $\sigma_j=\left<l_j,l_j+1\right>$ for $0\leq j\leq s$.  

\begin{ex} Hirzebruch Surfaces

In order to demonstrate this algorithm we will show that the minimal desingularization of ${\mathbb P}(1,1,n)$ (for $n\geq2$) is the $n^{th}$ Hirzebruch surface, ${\mathbb F}_n$.  ${\mathbb P}(1,1,n)$ is the complete toric surface given by $\Delta(1)=\left\{u_0,u_1,u_2\right\}$ where for $i=1,2$ $u_i=e_i$ constitute the standard basis for ${\mathbb N}={\mathbb Z}^2$ and $u_0=-e_1-ne_2$.  There are three maximal cones, $\sigma_0,\sigma_1,\sigma_2$ where $\sigma_i$ is spanned by $\Delta(1)\setminus u_i$.  Checking the determinants of each cone reveals \[ \begin{vmatrix}1& 0\\ 0& 1\end{vmatrix}=\begin{vmatrix}0& -1\\ 1& -n\end{vmatrix}=1\qquad\begin{vmatrix}-1& 1\\ -n& 0\end{vmatrix}=n\] so that the only singular cone is $\sigma_2=\left<u_0,u_1\right>$. 

In order to desingularize $\sigma_2$ take $v_1=u_1=e_1$ and $v_2=u_0=-e_1-ne_2$ and let $n_1=ae_1+be_2$.  Since $\det|v_1\: n_1|=\pm1$ we see that $b=\pm1$, so $n_1=ae_1\pm e_2$.  Now we need to find relatively prime integers $0\leq p_0<q_0$ such that \[ \begin{bmatrix}-1\\ -n\end{bmatrix}=p_0\begin{bmatrix}1\\ 0\end{bmatrix}+q_0\begin{bmatrix}a\\ \pm1\end{bmatrix}. \] From the bottom row we know that $q_0=n$ and $b=-1$.  The top row then yields $p_0=-1-an$.  The only way to satisfy $0\leq p_0<q_0$ is to have $a=-1$ and $p_0=n-1$.  Subsequently, $n_1=-e_1-e_2$.  

Next we must consider $\beta_0=\frac{q_0}{q_0-p_0}=\frac{n}{n-(n-1)}=n$.  Note that the Hirzebruch-Jung expansion for this is trivial, $b_1=n$.  Carrying out the rest of the algorithm then yields $l_1= l_0+n_1=e_1-e_1-e_2=-e_2$ as well as \[ n_2=(b_1-2)l_0+(b_1-1)n_1=(n-2)\begin{bmatrix}1\\ 0\end{bmatrix}+(n-1)\begin{bmatrix}-1\\ -1\end{bmatrix}=\begin{bmatrix}-1\\ -(n-1)\end{bmatrix}.\]  This leaves $l_2=l_1+n_2=-e_1-ne_2=v_2$, ending the algorithm.  

Compiling all of this data shows that the minimal desingularization of of ${\mathbb P}(1,1,n)$ is the complete toric surface whose fan in ${\mathbb N}={\mathbb Z}^2$ is given by \[ \Delta(1)=\left\{\begin{bmatrix}1\\ 0\end{bmatrix},\begin{bmatrix}0\\ 1\end{bmatrix},\begin{bmatrix}-1\\ -n\end{bmatrix},\begin{bmatrix}0\\ -1\end{bmatrix}\right\}\] which is easily recognized as the $n^{th}$ Hirzebruch surface, ${\mathbb F}_n$.  
\end{ex} 

\section{The Hirzebruch-Jung Continued Fraction}
\label{sec:cfrac}

In this section we break down the Hirzebruch-Jung continued fraction decompisition.  In particular we construct a method of determining each $b_j$ inductively from $\frac{q_0}{q_0-p_0}$.  

In order to do this we will construct a sequence of rational numbers, $\beta_0,\ldots\beta_{s-1}$, where each $\beta_j>1$.  These numbers will be defined by letting $\beta_j=[[b_{j+1},\ldots,b_s]]$.  Thus $\beta_0=\frac{q_0}{q_0-p_0}$ and $\beta_{s-1}=b_s$ is an integer.  It is important to note that this is the only integer in this sequence because once an integral $\beta_j$ is obtained, it has a trivial continued fraction of $\beta_j=[[b_{j+1}]]=b_{j+1}$, ending the process at $j=s-1$.  

\begin{lem}
\label{lem:inductbeta}

The fractions $\beta_j$ satisfy \[ \beta_{j+1}=\frac{1}{b_{j+1}-\beta_j}. \] 
\end{lem} 
\begin{proof} 
The proof is to simply write out the continued fraction for $\beta_j$ from the definition as \[ \beta_j=[[b_{j+1},\ldots,b_s]]=b_{j+1}-\frac{1}{[[b_{j+2},\ldots,b_s]]}=b_{j+1}-\frac{1}{\beta_{j+1}} \] and solve for $\beta_{j+1}$.  
\end{proof}

The next step of the process is to find a method to determine the values $b_j$ from the previous $\beta_j$'s.  

\begin{thm}
\label{thm:b=ceilbeta}

The integers from the Hirzebruch-Jung continued fraction satisfy $b_{j+1}=\left\lceil\beta_j\right\rceil$.  
\end{thm} 
\begin{proof} 
The proof is by induction, starting at $j=s-1$ and then decreasing $j$.  In the initial case the definition of the $\beta_j$'s says that $\beta_{s-1}=b_s$.  Subsequently $\beta_{s-1}$ is an integer.  It is then equal to its ceiling, proving the first case.  

For $j<s-1$, however, we know that $\beta_j$ is not an integer.  In this case, we assume the hypothesis that $b_{j+1}=\left\lceil\beta_j\right\rceil$ and try to prove the next step: $b_j=\left\lceil\beta_{j-1}\right\rceil$.  

Since $b_{j+1}=\left\lceil\beta_j\right\rceil$, and $b_{j+1}\geq2$ we see that $\beta_j>1$.  Subsequently $0<\frac{1}{\beta_j}<1$.  

Lemma \ref{lem:inductbeta} shows that $\beta_{j-1}=b_j-\frac{1}{\beta_j}$.  Solving for $\frac{1}{\beta_j}$ yields $b_j-\beta_{j-1}$.  Consequently $0< b_j-\beta_{j-1}<1$ so that $b_j=\left\lceil\beta_{j-1}\right\rceil$, and the theorem is proven.  
\end{proof}

This result will allow one to determine all of the $b_j$'s from $\beta_0$.  Substituting this result into Lemma \ref{lem:inductbeta} shows that \begin{equation}\label{eqn:inductbeta} \beta_{j+1}=\frac{1}{\left\lceil\beta_j\right\rceil-\beta_j}\end{equation} which allows one to construct the $\beta_j$'s using only the value of $\beta_0$.  Theorem \ref{thm:b=ceilbeta} lets one obtain the values in the continued fraction decomposition from this data.  

\begin{ex}$\beta_0=\frac{182}{27}$

In order to compute the continued fraction expanson of $\frac{182}{27}$ begin by noting that $b_1=\left\lceil\frac{182}{27}\right\rceil=7$.  Then $\beta_1=1/(7-\frac{182}{27})=\frac{27}{7}$, so that $b_2=\left\lceil\frac{27}{7}\right\rceil=4$.  Subsequently $\beta_2=1/(4-\frac{27}{7})=7$.  Since this is an integer, the sequence stops with $b_3=7$, so that $\frac{182}{27}=[[7,4,7]]$.  
\end{ex}

\section{Ending The Algorithm}
\label{sec:end} 

It turns out that the numbers $\beta_j$ are more than a simple intermediary sequence used to find $b_j$.  Consider what would happen if one were to start desingularizing the cone $\sigma$, then stop after obtaining some $l_j$ where $j\leq s$.  If one were to resume the process by desingularizing the new cone $\sigma'=\left<l_j,v_2\right>$ the resulting subdivisions must be the same as if the process had never been stopped and restarted.  

Not only must the $l_j$'s be the same, but the self-interesction numbers of the $T$-equivariant divisors $D(l_j)^2=-b_j$ must be the same as well.  This helps to give us a geometric interpretation of $\beta_j=[[b_{j+1},\ldots,b_s]]$.  Namely, if we were to attempt desingularizing the cone $\left<l_j, v_2\right>$, by finding appropriate $n_{j+1}\in{\mathbb Z}^2$ and $p_j,q_j\in{\mathbb Z}$, then $\beta_j$ would be $\frac{q_j}{q_j-p_j}$.  

There is, however, one exception to this: the non-singular cone.  Since the lattice points $l_0=v_1,l_1,\ldots,l_s,l_{s+1}=v_2$ give the minimal desingularization, attpempting to desingularize a non-singular cone would mean desingularizing $\left<l_s,v_2\right>$.  However, the sequence $\beta_j$ ends at $s-1$.  An attepmt to find $\beta_s$ from Lemma \ref{lem:inductbeta} would yield \[ \beta_s=\frac{1}{b_s-\beta_{s-1}}=\frac{1}{b_s-b_s} \] which is undefined.  

Let us now compare this to what happens when we attempt to desingularize this cone by finding $p_s$ and $q_s$.  The process begins by finding $n_{s+1}$ and $0\leq p_s<q_s$ such that $\det|l_s\; n_{s+1}|=\pm1$ and $l_{s+1}=p_sl_{s}+q_sn_{s+1}$.  Since $\left<l_s, l_{s+1}\right>$ is non-singular, taking $n_{s+1}=l_{s+1}$ will sasisfy the first requirement, while choosing $p_s=0$ and $q_s=1$ yields the second condition.  

Unfortunately, this does not match up with the rest of the algorithm.  In particular, it yields a value of $\beta_s=\frac{1}{1-0}=1$ while the algorithm yields and undefined $\beta_s$.  Furthermore, the algorithm itself says that $l_{j+1}=l_j+n_{j+1}$, while this choice of $p_s$, $q_s$ gives $l_{s+1}=n_{s+1}$ instead.  In particular this yields a value of $b_1=1$, which should correspond to the self-intersection number of the $T$-Weil divisor associated to the first new subdivision made.  Since the cone is already smooth there are no new subdivisions, so $b_1$ shouldn't even exist in this case.  

The way to address this problem is by changhing the requirement that $0\leq p_j<q_j$.  By considering $0<p_j\leq q_j$ instead, we can make the process of desingularizing an already non-singular cone consistent with the end of the algorithm.  

By changing this requirement, now when we try to find $p_s$, $q_s$ and $n_{s+1}$, we will find that $n_{s+1}=l_{s+1}-l_s$ satisfies \[ \det|l_s\; n_{s+1}|=\det|l_s\; l_{s+1}-l_s|=\det|l_s\; l_{s+1}|-\det|l_s\; l_s|=1 \] so that $l_s$ and $n_{s+1}$ form a ${\mathbb Z}$-basis for ${\mathbb N}$ as required.  In addition, the requirement $l_{s+1}=p_sl_s+q_sn_{s+1}$ results in $p_s=q_s=1$.  This corresponds to a value for $\beta_s$ of $\frac{q_s}{q_s-p_s}=\frac{1}{1-1}$ which is undefined, matching our result from the algorithm.  It is important to note that the condition $l_{s+1}=l_s+n_{s+1}$ is now satisfied as well.  

Thus, in the original algorithm, the natural way to describe the non-singular case is not the original formulation where $p_0=0$.  Rather, it corresponds to $p_0=q_0=1$.  Thus, by using a modification of the desingularization algorithm where $0<p_0\leq q_0$, we maintain consistency between ending the algorithm, and attempting to desingularize cones that are already smooth.  

\section{Desingularization of Weighted Projective Planes}
\label{sec:P(1,m,n)}

The goal of this section is to apply the breakdown of the Hirzebruch-Jung continued fraction in section \ref{sec:cfrac} to studying weighted projective planes.  Recall that a weighted projective plane, ${\mathbb P}(l,m,n)$ is the complete toric variety whose fan in ${\mathbb N}={\mathbb Z}^2$ is given by $\Delta(1)=\left\{ u_0,u_1,u_2\right\}$ where $lu_0+mu_1+nu_2=0$ and $\Delta(1)$ spans ${\mathbb N}$ as a ${\mathbb Z}$-module with standard basis $e_1,e_2$.  In the case where $l=1$ we can take $u_i=e_i$ for $i=1,2$ and $u_0=-me_1-ne_2$.  From here the lattice isomorphism given by $e_1\mapsto e_2$ and $e_2\mapsto e_1$ allows one to reduce to the case where $m<n$.  Also, the fact that $X(\Delta)$ is determined by $\sigma\otimes{\mathbb R}$ simplifies the problem by only having to consider the case where $m$, $n$ are relatively prime.  Subsequently, we may write $n=mk+r$ where $0<r<m$ and $m$, $r$ are relatively prime.  

This complete toric variety has three maximal cones, $\sigma_0$, $\sigma_1$, and $\sigma_2$ where each $\sigma_i$ is generated (over ${\mathbb R}_+$) by $\Delta(1)\setminus\{u_i\}$.  Since $\det|e_1\; e_2|=1$, we see that $U_{\sigma_0}$ is smooth.  On the other hand, since $\det|u_2\; u_0|=m$ and $\det|u_0\; u_1|=n$ we see that both $U_{\sigma_1}$ and $U_{\sigma_2}$ will be singular.  Since ${\mathbb P}(1,m,n)$ is normal, these singularities will have a codimension of at least 2.  Subsequently they are isolated points, fixed by the torus action on ${\mathbb P}(1,m,n)$ (i.e. the origin of each affine variety).  In terms of the classical homogeneous coordinates on ${\mathbb P}(1,m,n)$ these are the points $[0,1,0]$ and $[0,0,1]$ respectively.  

We will denote the minimal desingularization of ${\mathbb P}(1,m,n)$ by ${\mathbb D}(1,m,n)$, and its corresponding fan $\bar{\Delta}$.  We can begin to describe this surface by attempting to desingularize $\sigma_1=\left<u_2, u_0\right>$.  The following Lemma gets us as far down this road as we can without knowing more about $n$ and $m$.  

\begin{lem}
\label{lem:D(sigma1)}
When desingularizing $\sigma_1$ by taking $l_0=v_1=u_2$ and $l_{s+1}=v_2=u_0$, the first step has $l_1=-e_1-ke_2$.  
\end{lem}
\begin{proof}
Let $n_1=ae_1+be_2$.  Then the requirement that $\det|u_2\; n_1|=\pm1$ implies $a=\pm1$.  Now we need $0<p_0\leq q_0$ so that $u_0=p_0u_2+q_0n_1$.  By writing this out as \[ \begin{bmatrix}-m\\ -n\end{bmatrix}=p_0\begin{bmatrix}0\\ 1\end{bmatrix}+q_0\begin{bmatrix}\pm1\\ b\end{bmatrix} \] we see from the top row that $a=-1$ and $q_0=m$.  The bottom row then yields $p_0=-n-mb=-mk-r-mb$.  So in order to have $0<p_0\leq m$ we must take $b=-(k+1)$ so that $p_0=m-r$.  

Now that we have $\beta_0=\frac{q_0}{q_0-p_0}=\frac{m}{r}$ we see that $b_1=\left\lceil\frac{m}{r}\right\rceil$.  We also have $n_1=-e_1-(k+1)e_2$.  Plugging the last into the algorithm yields \[ l_1=l_0+n_1=\begin{bmatrix}0\\ 1\end{bmatrix}+\begin{bmatrix}-1\\ -(k+1)\end{bmatrix}=\begin{bmatrix}-1\\ -k\end{bmatrix}\] concluding the proof.  
\end{proof}

In order to address $\sigma_2=\left<u_0, u_1\right>$ we will take $l_0=v_1=u_1$ and $l_{s+1}=v_2=u_0$.  The following lemma tells us where the first two subdivisions occur.  

\begin{lem}
\label{lem:D(sigma2)}
When desingularizing $\sigma_2$ with $v_1=u_1$, the first two subdivisions occur at $l_1=-e_2$ and $l_2=-e_1-(k+1)e_2$.  
\end{lem}
\begin{proof}
We follow the approach of Lemma \ref{lem:D(sigma1)} by letting $n_1=ae_1+be_2$.  Since it must satisfy $\det|u_1\; n_1|=\det|e_1\; n_1|=\pm1$ we see that $b=\pm1$.  

Next we must find $0<p_0\leq q_0$ such that $v_2=p_0v_1+q_0n_1$, which we can write out as \[ \begin{bmatrix}-m\\ -n\end{bmatrix}=p_0\begin{bmatrix}1\\ 0\end{bmatrix}+q_0\begin{bmatrix}a\\ \pm1\end{bmatrix}. \] The bottom row reveals $b=-1$ and $q_0=n$, while solving the top row gives $p_0=-m-an$.  In order to have $p_0<q_0=n$ we must take $a=-1$, which gives $p_0=n-m$.  It also tells us that $n_1=-e_1-e_2$.  

Now that we have $\beta_0=\frac{q_0}{q_0-p_0}=\frac{n}{m}$, we use the fact that $n=mk+r$ to find $b_1=\left\lceil\beta_0\right\rceil=\left\lceil\frac{mk+r}{m}\right\rceil=k+1$.  Now we can proceed with the rest of the algorithm.  

First we find $l_1=l_0+n_1=e_1+(-e_1-e_2)=-e_2$.  Now the algorithm also tells us $n_2=(b_1-2)l_0+(b_1-1)n_1$.  Thus \[ n_2=(k-1)\begin{bmatrix}1\\ 0\end{bmatrix}+k\begin{bmatrix}-1\\ -1\end{bmatrix}=\begin{bmatrix}-1\\ -k\end{bmatrix}\] which can be plugged into $l_2=l_1+n_2$ to see $l_2=-e_1-(k+1)e_2$, proving the lemma.  
\end{proof}

It is well known that any smooth toric surface can be obtained from either ${\mathbb P}^2$ or a Hirzebruch surface ${\mathbb F}_n$ via a series of T-equivariant blowups \cite{Fult}.  The inclusion of $-e_2$ in Lemma \ref{lem:D(sigma2)} gives us an explicit demonstration of this fact for the smooth surface ${\mathbb D}(1,m,n)$.  This is because it is the complete toric variety given by \[ \bar{\Delta}(1)=\left\{\begin{bmatrix}0\\ 1\end{bmatrix},\begin{bmatrix}1\\ 0\end{bmatrix},\begin{bmatrix}0\\ -1\end{bmatrix},\begin{bmatrix}-1\\ -(k+1)\end{bmatrix},\ldots,u_j,\ldots\right\}\] where the $u_j$'s are all contained in the cone $\left<\begin{bmatrix}-1\\ -k\end{bmatrix}, \begin{bmatrix}-1\\ -(k+1)\end{bmatrix}\right>$.  Since this is a subcone of $\left<e_2, -e_1-(k+1)e_2\right>$, which is a single cone of ${\mathbb F}_{k+1}$ we see that the identity map on ${\mathbb N}={\mathbb Z}^2$ will give a map of toric varieties $\psi_{\left\lceil\cdot\right\rceil}:{\mathbb D}(1,m,n)\rightarrow{\mathbb F}_{\left\lceil\frac{n}{m}\right\rceil}$.  

Moreover, since $\bar{\Delta}(1)$ also contains $-e_1-ke_2$ and the $u_j$'s are contained in $\left<-e_1-ke_2, -e_2\right>$, we get a map $\psi_{\left\lfloor\cdot\right\rfloor}:{\mathbb D}(1,m,n)\rightarrow{\mathbb F}_{\left\lfloor\frac{n}{m}\right\rfloor}$ that is also induced by the the identity map on ${\mathbb N}$.  This shows that projection from a smooth toric surface to a Hirzebruch surface is not unique.  In fact, not even the Hirzebruch surface you map to is unique.  

While there are two different maps, it turns out that $\psi_{\left\lceil\cdot\right\rceil}$ has an added feature that $\psi_{\left\lfloor\cdot\right\rfloor}$ does not.  When desingularizing $\sigma_2$ we found the first term of the continued fraction to be $b_1=\left\lceil\frac{n}{m}\right\rceil=k+1$.  Because of this we know that the self-intersection number of the $T$-equivariant divisor associated to $l_1=-e_2$ is $\left[D(l_1)\right]^2=-b_1=-(k+1)$.  It can be seen in \cite{Faunt} and \cite{Hart} that the self-intersection number of the corresponding divisor on ${\mathbb F}_n$ is $-n$.  Thus, the map $\psi_{\left\lceil\cdot\right\rceil}:{\mathbb D}(1,m,n)\rightarrow{\mathbb F}_{k+1}$ preserves this quantity while $\psi_{\left\lfloor\cdot\right\rfloor}$ reduces it by one.  

\section{Applications}
\label{sec:applications}

Now that we know a few of the subsdivisions needed to desingularize ${\mathbb P}(1,m,n)$ in general, we will demonstrate the power of the techniques of Section \ref{sec:cfrac} by computing $\bar{\Delta}(1)$ for ${\mathbb D}(1,m,n)$ for a few explicit cases of $m$ and $n$.  In case $m=1$ we have seen that ${\mathbb D}(1,1,n)={\mathbb F}_n$ is a Hirzebruch surface, so we may restrict our attention to cases where $m\geq2$.  And in fact, the case $m=2$ is our first example.  

\begin{thm}
\label{thm:m=2}
${\mathbb D}(1,2,2k+1)$ is the complete toric surface whose fan in ${\mathbb Z}^2$ is given by \[ \Delta(1)=\left\{\begin{bmatrix}1\\ 0\end{bmatrix},\begin{bmatrix}0\\ 1\end{bmatrix},\begin{bmatrix}-1\\ -k\end{bmatrix},\begin{bmatrix}-2\\ -(2k+1)\end{bmatrix},\begin{bmatrix}-1\\ -(k+1)\end{bmatrix},\begin{bmatrix}0\\ -1\end{bmatrix}\right\}.\] 
\end{thm}
\begin{proof}
The easy way is to first use Lemmae \ref{lem:D(sigma1)} and \ref{lem:D(sigma2)} to see that these are the minimum lattice points necessary for any ${\mathbb D}(1,m,n)$.  Then simply check determinant of each maximal cone to check that it's smooth.  

It is a more straightforward proof (although a little more computationally intensive) to use the algorithm to explicitly desingularize the cones \begin{equation*}\begin{gathered} \sigma_1=\left<\begin{bmatrix}0\\ 1\end{bmatrix}, \begin{bmatrix}-2\\ -(2k+1)\end{bmatrix}\right>\\ \sigma_2=\left<\begin{bmatrix}-2\\ -(2k+1)\end{bmatrix},\begin{bmatrix}1\\ 0\end{bmatrix}\right>.\end{gathered}\end{equation*}
\end{proof}

This is the only case where the data from both Lemmae \ref{lem:D(sigma1)} and \ref{lem:D(sigma2)} complete the desingularization of ${\mathbb P}(1,m,n)$.  For larger values of $m$ at least one of the cones $\sigma_1$ or $\sigma_2$ will require further work.  We now turn our attention to the two cases where one the cones is desingularized by this data, but the other is not.  

The first case is when $m\geq2$ and $n=mk+1$.  In this case, $-e_1-ke_2$ will suffice to desingularize $\sigma_1$, while $\sigma_2$ will need to be blown-up $m$ times.  

\begin{thm}
\label{thm:sig1,r=1}
The only subdivision needed to desingularize $\sigma_1$ of ${\mathbb P}(1,m,n)$ when $n=mk+1$ is through the lattice point $-e_1-ke_2$.  
\end{thm}
\begin{proof}
Recall from the proof of Lemma \ref{lem:D(sigma1)} that by taking $l_0=e_2$ and $l_{s+1}=-me_1-(mk+1)e_2$ that $n_1=-e_1-(k+1)e_2$ and $\beta_0=\frac{m}{r}=m$.  Because this is an integer, the continued fraction is trivial and we see that $s=1$.  This means $l_{s+1}=l_2=-me_1-(mk+1)e_2$.  The algorithm gives the only remaining vector as $l_1=l_0+n_1=-e_1-ke_2$.  
\end{proof}

It is a straightforward enough task to verify that this is non-singular by checking that the determinant of each maximal cone is $\pm1$.  In fact, just as was the case when $m=2$ this would suffice to prove Theorem \ref{thm:sig1,r=1}.  In order to handle the other singular cone, $\sigma_2$ we must first determine the appropriate Hirzebruch-Jung decomposition.  Recall from the proof of Lemma \ref{lem:D(sigma2)} that for this cone the value $\beta_0=\frac{n}{m}$ was found.  

\begin{lem}
\label{lem:sig2,r=1,beta}
When $n=mk+1$, the Hirezebruch-Jung continued fraction, $[[b_1,\ldots,b_s]]$, of $\beta_0=\frac{n}{m}$ is $\beta_0=[[k+1,2,\ldots,2]]$ where $s=m$.  
\end{lem}
\begin{proof}
Since $\beta_0=\frac{n}{m}=\frac{mk+1}{m}$, Theorem \ref{thm:b=ceilbeta} is used to find $b_1=\left\lceil k+\frac{1}{m}\right\rceil=k+1$.  Subsquently, Lemma \ref{lem:inductbeta} tells us that $\beta_1=1/(b_1-\beta_0)=1/(k+1-\frac{n}{m})=\frac{m}{m-1}$.  

In order to complete the proof we will show that the remaining values of this sequence satisfy $\beta_j=\frac{m-(j-1)}{m-j}$.  In particular, this will mean that $b_{j+1}=\left\lceil\beta_j\right\rceil=2$ for $1\leq j\leq m-1$.  We proceed with a proof by induction, noting that it holds for $j=1$.  Suppose true for $j-1$.  Then \[ \beta_j=\frac{1}{\left\lceil\beta_{j-1}\right\rceil-\beta_{j-1}}=\frac{1}{2-\frac{m-(j-2)}{m-(j-1)}}=\frac{m-(j-1)}{m-j}\] as claimed.  

To see that $s=m$, note that \[ \beta_{m-1}=\frac{m-(m-2)}{m-(m-1)}=2\] is an integer.  The only integer in the rational sequence is the last one, $\beta_{s-1}$, so that $s-1=m-1$ concluding the proof.  
\end{proof}

Now that we have the continued fraction for $\beta_0=\frac{n}{m}$ in hand, we can continue with the desingularization of $\sigma_2$.  

\begin{thm}
\label{thm:sig2,r=1)}
When desingularizing $\sigma_2$ of ${\mathbb P}(1,m,mk+1)$ by taking $l_0=e_1$, and $l_{s+1}=-me_1-ne_2$, the subdivisions occur on the lattice points \[ l_j=\begin{bmatrix}-(j-1)\\ -\left[(j-1)k+1\right]\end{bmatrix}\] for $1\leq j\leq m+1=s+1$.  
\end{thm}
\begin{proof}
The proof is by induction.  Note that Lemma \ref{lem:D(sigma2)} shows this result is true for $j=1$ and $j=2$.  In that proof we already determined both $n_1=-e_1-e_2$ and $n_2=-e_1-ke_2$.  The algorithm tells us that the remaining $n_j$'s are determined by $n_{j+1}=(b_j-2)l_{j-1}+(b_j-1)n_j$.  However, Lemma \ref{lem:sig2,r=1,beta} says that $b_j=2$ for $j\geq2$, so that $n_{j+1}=(2-2)l_{j-1}+(2-1)n_j=n_j$.  Subsequently, for $j\geq2$ we have $n_j=n_2=-e_1-ke_2$.  

Now suppose that $l_j=-(j-1)e_1-\left[(j-1)k+1\right]e_2$.  Then \[ l_{j+1}=l_j+n_{j+1}=\begin{bmatrix}-(j-1)\\ -\left[(j-1)k+1\right]\end{bmatrix}+\begin{bmatrix}-1\\ -k\end{bmatrix}=\begin{bmatrix}-j\\ -(jk+1)\end{bmatrix} \] proving the theorem.  
\end{proof}

This completes the desingularization of ${\mathbb P}(1,m,mk+1)$.  In this case, Lemma \ref{lem:D(sigma1)} was sufficient in desingularizing $\sigma_1$, but Lemma \ref{lem:D(sigma2)} failed to completely desingularize $\sigma_2$.  We now turn our attention to the reverse case ($\sigma_2$ desingularizes quickly, but $\sigma_1$ does not).  Namely, we consider the surface ${\mathbb P}(1,m,mk+m-1)$.  

\begin{thm}
\label{thm:sig2,r=-1}
When desingularizing $\sigma_2$ of ${\mathbb P}(1,m,mk+m-1)$ the only necassary subdivisions are through $l_1=-e_2$ and $l_2=-e_1-(k+1)e_2$. 
\end{thm}
\begin{proof}
Recall from the proof of Lemma \ref{lem:D(sigma2)} that $\beta_0=\frac{n}{m}=\frac{mk+m-1}{m}$ so $b_1=\left\lceil\beta_0\right\rceil=k+1$.  Subsequently $\beta_1=\frac{1}{b_1-\beta_0}=m$ means $b_2=m$ is the last term of the decomposition.  In otherwords, $\beta_0=[[k+1,m]]$.  

The same proof also tells us that $n_1=-e_1-e_2$, so $l_1=-e_2$.  Then $n_2=(k-1)l_0+kn_1=-e_1-ke_1$ gives, $l_2=-e_1-(k+1)e_2$.  Subsequnently, \[ n_3=(m-2)l_1+(m-1)n_2=\begin{bmatrix} -(m-1)\\ -\left[(m-1)k+m-2\right]\end{bmatrix}\] so $l_{s+1}=l_3=l_2+n_3=-me_1-ne_2$ concluding the proof.  
\end{proof}

All that remains is to desingularize the cone \[\sigma_1=\left<\begin{bmatrix}0\\ 1\end{bmatrix}, \begin{bmatrix}-m\\ -(mk+m-1)\end{bmatrix}\right>.\]  In order to do this we will need the continued fraction of $\beta_0$ for this cone.  Recall that in the proof for Lemma \ref{lem:D(sigma1)} that this is $\beta_0=\frac{m}{r}=\frac{m}{m-1}$.  

\begin{lem}
\label{lem:sig1,r=-1,beta}
The Hirzebruch-Jung decomposition of $\beta_0=\frac{m}{m-1}$ is $b_j=2$ for all $1\leq j\leq s=m-1$.  
\end{lem}
\begin{proof}
We claim that the $\beta_j$'s are each $\frac{m-j}{m-(j+1)}$ for $0\leq j\leq m-1$.  Clearly it is true for $j=0$, and proceed by induction.  If $\beta_{j-1}=\frac{m-(j-1)}{m-j}$, then $b_j=\left\lceil\beta_{j-1}\right\rceil=2$.  Subsequently \[ \beta_j=\frac{1}{b_j-\beta_{j-1}}=\frac{1}{2-\frac{m-(j-1)}{m-j}}=\frac{m-j}{m-(j+1)}.\]  Note that this holds until $j=m-2$, when $\beta_{m-2}=\frac{m-(m-2)}{m-(m-1)}=2=b_{m-1}$ ends the sequence, so $s=m-1$.  
\end{proof}

With the necessary continued fraction in hand we may finish desingularizing ${\mathbb P}(1,m,mk+m-1)$ by desingularizing the cone $\sigma_1$.  

\begin{thm}
\label{thm:sig1,r=-1}
The minimal desingularization of the cone \\ $\sigma_1=\left<-e_2,-me_1-(mk+m-1)e_2\right>$ is given by  $\bar{\Delta}(1)=\left\{l_0,\ldots,l_m\right\}$ with \[ l_j=\begin{bmatrix}-j\\-(jk+j-1)\end{bmatrix}. \]  
\end{thm}
\begin{proof}
The proof of Lemma \ref{lem:D(sigma1)} shows that $n_1=-e_1-(k+1)e_2$.  However, since each $b_j=2$, we see that $n_{j+1}=(b_j-2)l_{j-1}+(b_j-1)n_j=n_j$.  Thus, $n_j=n_1=-e_1-(k+1)e_2$ for all $j$.  

We now proceed via induction, by noting that $l_0=e_2$ as needed.  Now suppose $l_{j-1}=-(j-1)e_1-\left[(j-1)k+j-2\right]e_2$.  Then \[ l_j=l_{j-1}+n_j=\begin{bmatrix}-(j-1)\\ -\left[(j-1)k+j-2\right]\end{bmatrix}+\begin{bmatrix}-1\\ -(k+1)\end{bmatrix}=\begin{bmatrix}-j\\ -(jk+j-1)\end{bmatrix}\] completing the proof.  
\end{proof}

We conclude the section by summarizing these results.  

\begin{cor}
\label{cor:D(1,m,mk+1)}
The minimal desingularization of ${\mathbb P}(1,m,mk+1)$ is the complete toric surface with \[ \Delta(1)=\left\{\begin{bmatrix}1\\ 0\end{bmatrix},\begin{bmatrix}0\\ 1\end{bmatrix},\begin{bmatrix}-1\\ -k\end{bmatrix},v_0,\ldots,v_m\right\}\] where $v_j=\begin{bmatrix}-j\\-(jk+1)\end{bmatrix}$.  
\end{cor}

\begin{cor}
\label{cor:D(1,m,mk+m-1)}
The minimal desingularization of ${\mathbb P}(1,m,mk+m-1)$ is the complete toric surface with \[ \Delta(1)=\left\{\begin{bmatrix}1\\ 0\end{bmatrix},\begin{bmatrix}0\\ -1\end{bmatrix},\begin{bmatrix}-1\\ -(k+1)\end{bmatrix},v_0,\ldots,v_m\right\}\] where $v_j=\begin{bmatrix}-j\\-(jk+j-1)\end{bmatrix}$.  
\end{cor}

\section{Programming the Hirzebruch-Jung Continued Fraction}
\label{sec:program}

The computation-minded reader will probably have already noted that as long as you are using something capable handling fractions (such as Maple), then equation \ref{eqn:inductbeta} and Theorem \ref{thm:b=ceilbeta} are sufficient for generating the continued fraction decomposition of $\beta_0$.  In other languages, however, there are problems with this approach.  In particular, the condition that the process stops when $\beta_{s-1}\in{\mathbb Z}$ becomes inapplicable because of rounding errors.  

The way around this is to construct an auxiliary sequence of integers $r_{-1},r_0,\ldots,r_{s-1}$ from which we may obtain each $\beta_j$.  Begin the sequence by letting $r_{-1}=q_0$.  Then generate each successive $r_j$ via the recursive relation \begin{equation}\label{eqn:rjdef} r_j=\frac{r_{j-1}}{\beta_j}.\end{equation}  The key here is to note that each $\beta_j$ can be obtained by taking the ratio of successive $r_j$'s, i.e. $\beta_j=\frac{r_{j-1}}{r_j}$.  It is not clear from this definition, however, that the $r_j$'s will be integers.  To see this we will need the following Theorem.  

\begin{thm}\label{thm:inductrj}
The sequence defined by Equation \ref{eqn:rjdef} satisfies \[ r_j=b_jr_{j-1}-r_{j-2}.\]  Subsequently, each $r_j$ is an integer.  
\end{thm}
\begin{proof}
First note that $\beta_j=\frac{r_{j-1}}{r_j}$.  On the other hand Lemma \ref{lem:inductbeta} says that $\beta_j=1/(b_j-\beta_{j-1})$.  Substituting the appropriate ratios of $r_j$'s for $\beta_j$ and $\beta_{j-1}$ then yields \[ \frac{r_{j-1}}{r_j}=\frac{1}{b_j-\frac{r_{j-2}}{r_{j-1}}}=\frac{r_{j-1}}{b_jr_{j-1}-r_{j-2}}.  \]  Setting the denominators equal proves the first statment of the theorem.  

To see that each $r_j$ is an integer note that both $r_{-1}=q_0$ and $r_0=\frac{r_{-1}}{\beta_0}=q_0-p_0$ are integers.  Then since each $b_j$ is an integer, each successive $r_j$ is integral and the theorem is proven.  
\end{proof}

In order implement this as a useful algorithm we must first eliminate the dependancy of this relation on the unknown values of $b_j$.  This can be done by using Lemma \ref{thm:b=ceilbeta} to substitue $\left\lceil\beta_{j-1}\right\rceil$ for $b_j$.  Then use the definition of $r_j$'s to replace $\beta_{j-1}$ with $\frac{r_{j-2}}{r_{j-1}}$.  Thus the $r_j$'s can be determined completely from $r_{-1}=q_0$ and $r_0=q_0-p_0$ by the recursive relation \begin{equation}\label{eqn:inductrj} r_j=\left\lceil\frac{r_{j-2}}{r_{j-1}}\right\rceil r_{j-1}-r_{j-2}.\end{equation}  

In order to develop a stop mechanism for this algorithm, recall that the fractional sequence $\beta_0,\ldots,\beta_{s-1}$ ended when $\beta_{s-1}$ was an integer.  Since $\beta_j=\frac{r_{j-1}}{r_j}$ this means that the integral sequence $r_{-1},\ldots,r_{s-1}$ ends when one $r_j$ divides the previous one.  

We conclude by giving the definition of a function written in the Python language for computing the Hirzebruch-Jung decomposition of a fraction $\frac{m}{n}>1$.  The input for the function {\tt HJcfrac} consists of the numerator and denominator of the fraction, given in any order.  There is no need to simplify the fraction (i.e., $m$ and $n$ can have a common factor).  The output of the function is then a list $b_1,\ldots,b_s$ of the integers in the continued fraction expansion of $\frac{m}{n}$.  

\begin{tabbing}
\noindent\tt >>> def \=\tt HJcfrac(m,n): \\
  \>\tt import math\\
  \>\tt r=[max(m,n),min(m,n)]\\
  \>\tt while \=\tt (r[-2]\% r[-1])>0:\\
  \>\>\tt  r.append((math.ceil(float(r[-2])/float(r[-1]))\\
  \>\>\tt  *r[-1])-r[-2])\\
  \>\tt for i in range(1,len(r):\\
  \>\>\tt  print ``b\_'',i,``='',int(math.ceil(\\
  \>\>\tt  float(r[j-1])/float(r[j])))
\end{tabbing}

\end{document}